\documentclass[12pt]{article}   
   
\usepackage{latexsym,amsfonts,amsmath,epsfig,tabularx,amsthm}   
   
\theoremstyle{definition}    
   
\newtheorem{thm}{Theorem}

\newcommand\reals{\mathbb{R}}   
\newcommand\R{{\mathbb{R}}}   
\renewcommand\natural{\mathbb{N}}

\newcommand\eps{\varepsilon}

\newcommand{\e}{\varepsilon}

\newcommand{\infp}{\operatornamewithlimits{inf\phantom{p}}}

\newcommand{\vol}{{\rm vol}}   
\newcommand{\conv}{{\rm conv}}   
   
\renewcommand{\P}{\mathbb P}   
\newcommand{\E}{\mathbb E}   
   
\newcommand\dint{{\rm d}}

\newlength{\fixboxwidth}   
\title{The Curse of Dimensionality for\\    
Monotone and Convex Functions of Many Variables}     

\date{September 15, 2010} 
 
\author{Aicke Hinrichs\footnote{This author  
was supported by the DFG Heisenberg grant HI 584/3-2.}, Erich Novak\footnote{This
author was partially supported by the DFG-Priority Program 1324}  \\  
Mathematisches Institut, Universit\"at Jena\\  
Ernst-Abbe-Platz 2, 07740 Jena, Germany\\  
email: a.hinrichs@uni-jena.de,  
erich.novak@uni-jena.de\\  
\qquad  
\\  
Henryk Wo\'zniakowski\footnote{This author was partially   
supported by the National Science  
Foundation.  
}\\  
Department of Computer Science, Columbia University,\\  
New York, NY 10027, USA, and\\  
Institute of Applied Mathematics, University of Warsaw\\  
ul. Banacha 2, 02-097 Warszawa, Poland\\  
email:\ henryk@cs.columbia.edu}  
 
\begin{document}   
   
\maketitle   
   
\begin{abstract}   
We study the integration and approximation problems for    
monotone and convex  
bounded functions that depend on $d$ variables,     
where $d$ can be arbitrarily large.    
We consider the worst case error    
for algorithms that use finitely many   
function values.    
We prove that these problems suffer from the curse of dimensionality.   
That is, one needs exponentially many (in $d$)    
function values to achieve an error $\e$.   
\end{abstract}    
   
\section{Introduction}\label{s2}   
Many multivariate problems suffer from the curse of dimensionality.  
A partial list of such problems can be found in e.g., \cite{NW08,NW10}.  
The phrase \emph{curse of dimensionality}   
was coined by Bellman already in 1957 and means   
that the complexity\footnote{By complexity we mean the minimal cost  
of computing an $\e$-approximation. The complexity is bounded from  
below by the information complexity which is defined as the minimal  
number of function values needed to compute an $\e$-approximation.  
In this paper we prove that even the information complexity suffers from  
the curse of dimensionality.} of a $d$-variate problem is an exponential  
function in $d$. This is usually proved for multivariate problems   
defined on the unit balls of normed linear spaces. We stress  
that the curse of dimensionality may hold independently of the   
smoothness of functions and may hold even for analytic functions.  
  
The choice of the unit ball as the domain of a multivariate problem  
is not essential and can be slightly generalized.   
What is important and heavily used in the proof is  
that the domain $F_d$ of the $d$ variate problem is balanced ($f\in  
F_d$ implies $-f\in F_d$) and convex ($f_1,f_2\in F_d$ and $t\in[0,1]$  
imply that $tf_1+(1-t)f_1\in F_d$).   
It is not clear if the curse of dimensionality may hold for  
domains $F_d$ being not balanced or not convex. 
   
In this paper we study classes of monotone and convex $d$-variate  
bounded functions. Such classes are obviously \emph{not} balanced and the   
previous analysis to prove the curse of dimensionality does not apply.   
We study the integration problem and the approximation problem   
in the $L_p$ norm with $p\in[1,\infty]$.  
We consider the worst case setting and  
algorithms that use finitely many function values. In particular, we  
ask what is the minimal number of $d$-variate function values 
that is needed to achieve an error $\e$.   
  
It turns out that the approximation problem in the $L_p$ norm for both  
monotone and convex functions is no easier than the integration  
problem. This means that lower error bounds for integration  
also hold for approximation. Hence, it is enough to prove the curse of  
dimensionality for the integration problem.   
  
The integration problem for monotone functions has been studied by  
Papageorgiou~\cite{P93}, and for convex functions by  
Katscher, Novak and Petras~\cite{KNP96}. They proved the optimal rate  
of convergence and provided lower and upper bounds on the $n$th  
minimal error. From these bounds we can conclude the lack of some  
tractability properties defined later,    
but cannot conclude whether the curse of dimensionality holds.     
  
In this paper we prove that for both monotone and convex functions,   
the curse of dimensionality holds for the   
integration problem and therefore also holds for the approximation  
problem in the $L_p$ norm.   
The proof relies on identifying  ``fooling''  
functions $f^{-}$ and $f^{+}$ which are both monotone or both convex, 
which share the same $n$ function   values used by an algorithm, and whose  
integrals differ as much as possible. Here ``as much as possible'' 
means that the error is at most $\e$ only if $n$ is exponentially  
large in $d$. The fooling functions for the monotone  
class take only values $0$ or $1$ depending   
on the points used by an algorithm.  
The fooling functions for the convex class are $f^{-}=0$ and $f^{+}$   
is chosen such that it vanishes at $n$ points used by an algorithm,  
and its integral is maximized. Using the results of  
Elekes~\cite{E86} and Dyer, F{\"u}redi and McDiarmid~ \cite{DFM90}    
on random volumes of cubes, we prove that the integral of $f^{+}$    
is of order $1$ for large $d$, if $n$ is smaller than, say, $(12/11)^d$.  
 
Restricting the algorithms for the integration problem to use  
only function values is quite natural.  
However, for the approximation problem it would be also interesting   
to consider algorithms that use finitely many arbitrary linear functionals.  
We believe that the $L_p$ approximation problem still suffers from the 
curse of dimensionality for this general information, and   
pose this question as an open problem. The paper by   
Gilewicz, Konovalov and Leviatan~\cite{GKL06}  
may be relevant in this case. This paper presents  
the order of convergence for the approximation problem   
for $s$-monotone functions (in one variable).  
  
We finally add a comment on the worst case setting used in this paper.  
Since integration for monotone and convex classes suffers from the  
curse of dimensionality in the worst case setting, it seems natural to  
switch to the randomized setting where algorithms can use function  
values at randomized sample points.   
Now we can use the classical Monte Carlo algorithm.    
Since all monotone and convex integrands are bounded by one,  
the error bound of Monte Carlo is $n^{-1/2}$, without any additional    
constant. Hence, $\e^{-2}$ function values at randomized sample points  
are enough to guarantee a randomized error $\e$. This means that   
the integration problem for both monotone and convex functions is    
\emph{strongly polynomially tractable}\footnote{This means 
  that~\eqref{poltra} holds with $q=0$. 
In this case we can choose $C=1$ and $p=2$.} 
 in the randomized setting.    
The exponent $2$ of $\e^{-1}$  is optimal since the    
optimal orders of convergence for randomized algorithms are   
$n^{-1/2-1/d}$  for monotone functions, see \cite{P93},    
and $n^{-1/2-2/d}$  for convex functions, see \cite{KNP96}.    
Hence, for large $d$ we cannot guarantee a randomized error $\e$ with   
$\e^{-p}$ function values with $p<2$.  
This proves that the switch for the worst case setting to the randomized  
setting breaks the curse of dimensionality for the integration problem   
defined for monotone and convex functions.   
   
Not much seems to be known about the $L_p$ approximation problem  
in the randomized setting  for monotone or convex functions.  
It is not clear if we still have the curse of dimensionality in   
the randomized setting. We pose this as another open problem.  
  
\section{Integration}  
We mainly study the integration problem, i.e., we want to approximate   
$$   
{\rm INT}_d(f) = \int_{[0,1]^d} f(x) \, {\rm d} x ,    
$$   
for bounded functions    
$f:[0,1]^d  \to [0,1] $ that are monotone   
(more precisely, non-decreasing in each variable $x_j$    
if the other variables are fixed)    
or convex. Hence, we consider the classes    
$$   
F_d^{\,\rm mon} = \{ f: [0,1]^d \to [0,1]   
\mid \ f \text{ is monotone} \}    
$$   
and    
$$   
F_d^{\, \rm con}   
= \{ f: [0,1]^d \to [0,1] \mid \ f \text{ is convex} \}.    
$$   
We approximate the integral ${\rm INT}_d(f)$ by algorithms   
$A_n$ that use    
information about $f$   
given by  
$n$ function values.    
Hence, $A_n$ has the form    
\begin{equation}   \label{form}   
A_n(f)=\phi_n(f(t_1), f(t_2), \dots, f(t_n)),    
\end{equation}    
where $n$ is a nonnegative integer, $\phi_n:\reals^n\to \R $ is an   
arbitrary mapping, and the choice of arbitrary     
sample points $t_j \in [0,1]^d$    
can be adaptive. That is, $t_j$ may depend on the   
already computed values    
$f(t_1), f(t_2), \dots, f(t_{j-1})$.   
For $n=0$, the mapping $A_n$ is a constant real number.   
More details can be found in e.g.,~\cite{No88,NW08,NW10,TWW88}.    
   
We define the $n$th minimal error   
of such approximations in the worst case setting as    
$$   
e_n^{\rm int} (F_d)  = \infp_{A_n} \ \sup_{f \in F_d}   
| {\rm INT}_d (f) - A_n (f) | \ \ \ \text{for }\ \ \   
F_d\in\{F^{\,\rm mon}_d,F^{\,\rm con} _d\}.   
$$   
For $n=0$, it is easy to see that the best    
algorithm is $A_0=\tfrac12$ for  the two classes considered    
in this paper, and we obtain   
$$   
e_0^{\rm int}(F_d^{\,\rm mon}) = e_0^{\rm int}(F_d^{\,\rm con}) = \tfrac12   
\ \ \ \text{for all}\ \ \  d \in \natural.    
$$   
   
Hence, the integration problems are well scaled     
and it is enough to study  the absolute error.    
The~\emph{information complexity} is the inverse function of   
$e_n^{\rm int}(F_d)$    
given by   
$$   
n^{\rm int}(F_d, \e) = \min \{ n \mid \ e_n^{\rm int}(F_d) \le \e \} \ \ \   
\text{for }\ \ \   
F_d\in\{F^{\,\rm mon}_d,F^{\,\rm con} _d\}.   
$$   
It is trivial that $n^{\rm int}(F_d^{\,\rm mon}, \e)=  
n^{\rm int}(F_d^{\,\rm con},\e)=0$    
for all $\e \ge \tfrac12$.    
   
\section{Known and new results}    
   
The integration problems for monotone and for    
convex  functions were studied    
before, we refer to the paper    
by Papageorgiou~\cite{P93} for monotone functions,    
and to the paper by   
Katscher, Novak and Petras~\cite{KNP96} for convex functions.   
Here we mention some of the known results and indicate    
our new results concerning the \emph{curse of dimensionality}.    
   
For the class $F_d^{\,\rm mon}$ of monotone functions   
it was proved by    
Papageorgiou~\cite{P93} that   
$$   
e_n^{\rm int}(F_d^{\,\rm mon}) = \Theta ( n^{-1/d}) .   
$$   
Hence, the \emph{optimal order of convergence}    
is $n^{-1/d}$. More precisely, it is proved in \cite{P93} that    
there are some positive numbers $c,C$  
independent of $n$ and $d$   
such that for all $d, n \in\natural$ we have     
\begin{equation}  \label{bound1}    
c\,d^{-1}\, n^{-1/d}\le    
e_n^{\rm int} (F_d^{\,\rm mon}) \le C\,{d} \, n^{-1/d}.   
\end{equation}    
  
It is interesting to note that the ratio    
between the upper and the lower bound is of the order    
$d^{\,2}$ , i.e.,  it is  polynomial in $d$,    
not exponential as it is   
the case for many other spaces.   
   
The bound~\eqref{bound1} yields   
$$   
\left\lceil \left(\frac{c}{d\,\e}\right)^d\right\rceil   
\le n^{\rm int}(F_d^{\,\rm mon},\e)\le   
\left\lceil \left(\frac{C\,d}{\e}\right)^d\right\rceil.    
$$   
{}From this we conclude that {\em polynomial tractability} and    
even {\em weak tractability} do not   
hold. That is, it is \emph{not} true that there are non-negative $C,q,p$   
such that for all $d\in \natural$ and $\e\in(0,\tfrac12)$ we have    
\begin{equation}\label{poltra} 
n^{\rm int}(F_d^{\,\rm mon},\e)\le C\,d^{\,q}\e^{-p}\ \ \    
\text{(polynomial tractability)},    
\end{equation} 
as well as it is \emph{not} true that   
$$   
\lim_{\e^{-1}+d\to\infty}\frac{\ln\, n^{\rm int}(F_d^{\,\rm   
    mon},\e)}{\e^{-1}+d}=0   
\ \ \ \text{(weak tractability)}.   
$$    
Nevertheless, the lower bound on $n^{\rm int}(F_d^{\,\rm mon},\e)$  
is useless     
for a fixed $\e >0$ and large $d$, since    
for $d\ge c/\e$ we do not obtain    
a bound better than $n^{\rm int}(F_d^{\,\rm mon}, \e) \ge 1$.   
Thus, it is not clear whether    
the information complexity    
$n^{\rm int}(F_d^{\,\rm mon}, \e)$ is exponential   
in $d$ for a fixed $\e\in (0,\tfrac12)$.    
In this paper we will prove that   
$$   
n^{\rm int}(F_d^{\,\rm mon}, \e)    
\ge 2^d \, (1-2\e)\ \ \ \text{for all}   
\ \ \ d\in\natural,\ \e\in (0,\tfrac12).    
$$   
This means that $n^{\rm int}(F_d^{\,\rm mon}, \e)$   
is indeed exponential in $d$,   
that is the integration problem   
suffers from the \emph{curse of dimensionality}.  
  
\vskip 1pc   
We now turn to the class $F_d^{\,\rm con}$ of convex functions.    
It was proved by Katscher, Novak and Petras~\cite{KNP96} that   
$$   
e_n^{\rm int} (F_d^{\,\rm con}) = \Theta ( n^{-2/d}) .   
$$   
Again, the \emph{optimal order of convergence}    
is known, now it is  $n^{-2/d}$. More precisely, it was proved    
in~\cite{KNP96} that there are some positive numbers $c_d,C$,   
with $c_d$ being exponentially small in   
$d$ whereas $C$ is independent of $d$, such that    
we have for all $n \in \natural$   
\begin{equation}\label{bound2}   
c_d \, n^{-2/d} \le e_n^{\rm int}(F_d^{\,\rm con}) \le C \, d \, n^{-2/d}.   
\end{equation}   
The bound~\eqref{bound2} yields    
$$   
\left\lceil \left(\frac{c_d}{\e}\right)^{d/2}\right\rceil   
\le n^{\rm int}(F_d^{\,\rm con},\e)\le   
\left\lceil \left(\frac{C\,d}{\e}\right)^{d/2}\right\rceil.    
$$   
{}From this we conclude that polynomial tractability does not hold.    
The lower bound in~\eqref{bound2} is useless     
for a fixed $\e >0$ and large $d$, and therefore it is not clear if we   
have weak tractability or the curse of dimensionality.   
In this paper we will prove that there exists $\e_0\in(0,1/4)$ such   
that    
$$   
n^{\rm int} (F_d^{\,\rm con},\e) 
\ge \frac1{2(d+1)}\left(\frac{11}{10}\right)^d   
\ \ \ \text{for all}\ \ \   
d\in\natural,\ \e\in(0,\e_0].   
$$   
Hence, the integration problem also suffers  from 
the curse of dimensionality for convex functions.    
   
\section{The class of monotone functions}    
  
We consider integration for monotone functions.    
Assume that $A_n$ is an arbitrary (possibly adaptive) algorithm for    
the class $F_d^{\,\rm mon}$. For $x=[x_1,x_2,\dots,x_d]\in [0,1]^d$,    
consider the ``fooling'' function   
$$   
f^*(x)=\begin{cases}   
0 & \text{if}\  \sum_{k=1}^d x_k < d/2,\\   
1 & \text{if}\  \sum_{k=1}^d x_k \ge d/2.   
\end{cases}   
$$    
Obviously, $f^* \in F_d^{\,\rm mon}$ and therefore the algorithm    
$A_n$ will use function values    
$$   
f^*(t_1), f^*(t_2), \dots , f^*(t_n)   
$$   
for some sample points $t_j\in[0,1]^d$.    
Since the algorithm $A_n$ can \emph{only}    
use the computed function values, we obtain    
$$   
A_n(f) = A_n(f^*)   
$$   
for all $f \in F^{\,\rm mon}_d$ if    
$f(t_k)=f^*(t_k)$ for $k=1,2, \dots , n$.    
   
Take first the case $n=1$.    
Suppose first that $f^*(t_1)=0$, i.e., $\sum_{j=1}^dt_{1,j}<d/2$    
for $t_1=[t_{1,1},t_{1,2},\dots,t_{1,d}]$. Define     
$f^-=0$ and the function    
$$   
f^+(x)=\begin{cases}   
0&  \text{if}\  x \le t_1\ \text{(in every coordinate)},\\    
1&  \text{otherwise.}   
\end{cases}   
$$   
Then $f^-,f^+\in F_d^{\,\rm mon}$ and they yield the same information    
as $f^*$, i.e.,    
$$   
f^-(t_1)=f^+(t_1)=f^*(t_1)=0.   
$$   
   
Using the standard proof technique it can be checked that   
$$   
\max_{y\in[0,1]^d,\ \sum_{j=1}^dy_j\le d/2}\prod_{j=1}^dy_j=   
\max_{y\in[0,1]^d,\ \sum_{j=1}^dy_j\ge d/2}\prod_{j=1}^d(1-y_j)=2^{-d}.   
$$   
Then   
$$   
{\rm INT}_d(f^+)   
=1-{\rm INT}_d(1-f^+)=1-\int_{x\le t_1}{\rm d}x=   
1-\prod_{j=1}^dt_{1,j}.   
$$    
This implies that     
\begin{equation}   \label{eq2}    
{\rm INT}_d(f^+)- {\rm INT}_d(f^-) \ge 1 - 2^{-d} .   
\end{equation}    
The case with $f^*(t_1)=1$ is similar.    
Now take $f^+=1$ and    
$$   
f^-(x) = \begin{cases}   
1 & \text{if}\ x \ge t_1,\\   
0 & \text{otherwise}.   
\end{cases}   
$$   
Again $f^+$ and $f^-$ are from $F_d^{\,\rm mon}$ and they    
yield the same information    
as $f^*$. We also obtain~\eqref{eq2}.    
We estimate the error of $A_1$ on the whole class    
$F_d^{\,\rm mon}$ by    
\begin{eqnarray*}   
\sup_{f\in F^{\,\rm mon}_d}|{\rm INT}_d(f)-A_n(f)|&\ge&   
\max\left(|{\rm INT}_d(f^+)-A_n(f^*)|,|{\rm   
    INT}_d(f^-)-A_n(f^*)|\right)\\   
&\ge&   
\tfrac12\left({\rm INT}_d(f^+)-{\rm INT}_d(f^-)|\right)   
\ge \tfrac12\left(1-2^{-d}\right).   
\end{eqnarray*}    
Since this holds for all algorithms, we conclude that   
$$   
e_1(F^{\,\rm mon}_d,\e)\ge \tfrac12\left(1-2^{-d}\right).   
$$   
   
The general case with $n \in \natural$ is similar.    
Assume that $\ell$ of the function values yield    
$f^*(t_k)=0$ while $n-\ell$ function values yield    
$f^*(t_k)=1$. Without loss of generality, we may assume that   
\begin{eqnarray*}   
f^*(t_j)&=&0\ \ \ \text{for}\ \ \ j=1,2,\dots,\ell,\\   
f^*(t_j)&=&1\ \ \ \text{for}\ \ \ j=\ell+1,\ell+2,\dots,n.   
\end{eqnarray*}   
Define the two functions,   
$$   
f^+(x)=\begin{cases}   
0& \text{if}\ x\le t_1\ \text{or}\ x\le t_2\ \text{or}\ \ldots   
\ \text{or}\ x\le t_\ell,\\   
1& \text{otherwise}.   
\end{cases}   
$$   
and    
$$   
f^-(x)=\begin{cases}   
1& \text{if}\ x\ge t_{\ell+1}\ \text{or}\   
x\ge t_{\ell+2}\ \text{or}\ \ldots   
\ \text{or}\ x\ge t_n,\\   
0& \text{otherwise}.   
\end{cases}   
$$   
Then $f^+, f^- \in F^{\,\rm mon}_d$ with    
$$   
f^+(t_k)=f^-(t_k)=f^*(t_k)\ \ \   
\text{for all}\ \ \ k=1,2, \dots , n.   
$$   
Furthermore, we have    
$$   
{\rm INT}_d(f^-) \le \sum_{j=1}^{n-\ell}\int_{x \ge t_{\ell+j}} 1 \,   
{\rm d}x\le (n-\ell)2^{-d}.    
$$   
Similarly it is easy to show that    
${\rm INT}_d(f^+) \ge 1- 2^{-d} \cdot \ell$,  so that   
$$   
{\rm INT}_d(f^+)-{\rm INT}_d(f^-) \ge 1 - 2^{-d} \cdot n .   
$$   
Therefore the worst case error of $A_n$ is at least    
$\tfrac12 (1-2^{-d} n)$. Since this holds for an arbitrary $A_n$    
we also have   
$$   
e_n(F_d^{\,\rm mon})\ge \tfrac12\left(1-2^{-d}n\right).   
$$   
This leads to the following theorem.     
\begin{thm}   
\label{thm:_tracmon}   
For each fixed $\e\in(0, \tfrac12)$, the information    
complexity is at least    
$$   
n^{\rm int}(F_d^{\,\rm mon}, \e) \ge 2^d \, (1-2\e)\ \ \ \text{for all}   
 \ \ \ d\in\natural.    
$$   
Thus, the integration problem for monotone functions    
suffers from the curse of dimensionality.    
\end{thm}    
   
\section{The class of convex functions}    
   
We now consider integration for convex function and    
prove the curse of dimensionality.  
  
\begin{thm}  
\label{thm:_tracconv}  
There exists $\eps_0\in(0,\tfrac12)$   
such that for each fixed $\eps\in(0,\eps_0)$   
the information complexity is at least  
$$  
n^{\rm int}  
(F_d^{\,\rm con}, \e) \ge \frac{1}{d+1} \left( \frac{11}{10} \right)^d   
\left( 1-\frac{\eps}{\eps_0} \right)     
\ \ \ \text{for all} \ \ \ d\in\natural.    
$$   
Thus, the integration problem of convex functions    
suffers from the curse of dimensionality.    
\end{thm}  
  
The idea of the proof is as follows. Assume again that we have an   
arbitrary (possibly adaptive) algorithm $A_n$ for the class   
$F_d^{\,\rm con}$.  
For the zero function $f^-=0$ the algorithm $A_n$ uses function   
values at certain sample points $x_1, x_2, \dots,x_n$.  
This implies that $A_n$ uses the same sample points $x_1,x_2, \dots,x_n$   
for any function $f$ from $F_d^{\,\rm con}$ with   
$$f(x_1)= f(x_2) = \dots=f(x_n)=0.  
$$  
In particular, let $f^+$ be the largest such   function,  
$$  
f^+(x)=\sup\{f(x)\,|\ \ f(x_j)=0,\ j=1,2,\dots,n,\ f\in F_d^{\,\rm  
  con}\,\}.  
$$  
Clearly, $f^+\in F_d^{\,\rm con}$, $f^+(x_j)=0$ for $j=1,2,\dots,n$,  
$f(x)\ge0$ for all $x\in[0,1]^d$, and $f^+$ has the  
maximal value of the integral among such functions.   
The integral ${\rm INT}_d(f^+)$ is the volume of the subset   
under the graph of the function $f^+$.  
This subset under the graph is the complement in $[0,1]^{d+1}$ of  
the convex hull of the points $(x_1,0), (x_2,0),  
\dots,(x_n,0)\in[0,1]^{d+1}$   
and $[0,1]^d \times \{1\} \subset [0,1]^{d+1}$.    
Denoting this convex hull by $C$, we obtain  
$$  
{\rm INT}_d(f^+) = 1 - \vol_{d+1} (C).  
$$  
  
Since the algorithm $A_n$ computes the same result for the   
functions $f^-$ and $f^+$ but ${\rm INT}_d(f^-) = 0$ we conclude that  
$A_n$ has error at least   
$$  
\tfrac{1}{2} \big( 1 - \vol_{d+1} (C) \big)  
$$  
on one of these functions.  
Theorem \ref{thm:_tracconv} now follows   
directly from the next theorem which gives an estimate of the   
volume of the set $C$ by setting $\eps_0=t_0/2$.   
  
\begin{thm}       
Let $P$ be an $n$-point set in $[0,1]^d \times \{0\}$.    
Then the $(d+1)$-dimensional volume of the convex hull $C$  
of $P \cup \big( [0,1]^d \times \{1\} \big)$ is at most    
$$   
\vol_{d+1} (C) \le (1-t_0) + (d+1)\, n\, t_0 \left( \frac{10}{11} \right)^d  
$$   
for some $t_0\in(0,1)$ independent of $d$ and $n$.    
\end{thm}   
   
\begin{proof}   
Let $Q=[0,1]^d$ and $Q_t=[0,1]^d\times \{t\}  
\subset \R^{d+1}$ for $t\in [0,1]$.    
Let $P\subset Q_0$ be an $n$-point set and let $C$ be the   
convex hull of $P \cup Q_1$.  We want to show that   
$$    
\vol_{d+1} (C)   \le (1-t_0) + (d+1)\, n\, t_0  
\left( \frac{10}{11} \right)^d.   
$$   
Let $C_t=C\cap Q_t$ be the slice of $C$ at height $t$.  
For a point $z=(z_1, z_2, \dots,z_d,z_{d+1}) \in\R^{d+1}$ let   
$\overline z =(z_1, z_2, \dots,z_d)$ be its projection onto   
the first $d$ coordinates.  
Similarly, for a set $M\subset \R^{d+1}$, let $\overline M$   
be the set of all points $\overline z$ with $z\in M$.   
  
Since   
$$    
\vol_{d+1} (C)  = \int_0^1 \vol_d(C_t) \, {\rm d} t =   
\int_0^1 \vol_d(\overline C_t) \, {\rm d} t  \le (1-t_0) + \int_0^{t_0}   
\vol_d(\overline C_t) \, {\rm d} t,   
$$   
it is enough to prove that  
$$   
\vol_d(\overline C_t) \le (d+1)\, n  \left( \frac{10}{11} \right)^d  
\ \ \ \mbox{for all}\ \ \ t\in[0, t_0].  
$$  
  
Carath\'eodory's theorem states that any point in the convex   
hull of a set $M$ in $\R^d$ is already   
contained in the convex hull of a subset of $M$ consisting   
of at most $d+1$ points.   
Hence, every point of $P$ is contained in the convex hull   
of $d+1$ vertices of $Q_0$.   
It follows that it is enough to show that   
\begin{equation}   
\label{eq:_1}   
\vol_d(\overline C_t) \le  n  \left( \frac{10}{11} \right)^d  
\end{equation}   
whenever $P$ is an $n$-point set of such vertices of $Q_0$.   
So we assume now that $P$ is such a set.   
  
Let   
$$  
w_t = ((1+t)/2, (1+t)/2, \dots,(1+t)/2,t) \in Q_t.  
$$   
For each vertex $v\in P$, let $B_v\subset Q_0$ be the intersection   
of the ball with center   
$\tfrac12(w_0 +v)$ and radius $\tfrac{1}{2} \Vert w_0 - v\Vert$   
with $Q_0$.   
Observe that $C_0$ is the convex hull of $P$.    
By Elekes' result from~\cite{E86},   
$$   
C_0 \subset \bigcup_{v \in P} B_v.  
$$   
It follows that   
$$   
C = \conv (P \cup Q_1)  \subset \bigcup_{v \in P} \conv( B_v \cup Q_1)   
$$   
since each point in this convex hull lies on a  segment   
between a point in some $B_v$ and a point in $Q_1$.   
Since all sets  $\conv( B_v \cup Q_1)$ are congruent,    
the inequality \eqref{eq:_1} immediately follows if we show that   
\begin{equation}   
\label{eq:_2}   
\vol_d ( \overline D_t ) \le \left( \frac{10}{11} \right)^d  
\ \ \ \mbox{for all}\ \ \ t\in[0,t_0],   
\end{equation}   
where $D_t= \conv( B_v \cup Q_1) \cap Q_t$   
is the section of the convex hull at height $t$.   
We can now restrict ourselves to the case that $v$ is a fixed    
vertex in $P$, say $v=(0, 0, \ldots,0,0)$.   
    
Let $O$ be the origin in $\R^d$.   
Let $E_t\subset Q$ be the intersection of the ball with    
center $\tfrac12{\overline w_t}$ and diameter   
$\Vert  \overline w_t \Vert$ with $Q$.   
Then $\overline D_t \subset E_t$, so \eqref{eq:_2} is proved once we show   
\begin{equation}   
\label{eq:_3}   
\vol_d ( E_t ) \le \left( \frac{10}{11} \right)^d  
\ \ \ \mbox{for all}\ \ \ t\in[0,t_0].     
\end{equation}    
  
To this end we follow the approach from \cite{DFM90}.   
Set $2s=\tfrac12(1+t)$. Then   
$$    
\vol_d ( E_t ) = \P \Big( \sum_{j=1}^d (X_j-s)^2 \le d s^2 \Big)    
$$   
where $X_1,X_2, \dots,X_d$ are independent uniformly distributed in $[0,1]$.   
We now use Markov's inequality   
$$  
\P ( |Y| \ge a ) \le \frac{\E (|Y|)}{a},  
$$  
which holds for all real random variables $Y$ and all $a> 0$.  
We take $a=1$ and   
$$  
Y =   \exp \Bigg( \alpha \Big( d s^2 - \sum_{j=1}^d (X_j-s)^2 \Big) \Bigg),  
$$   
and conclude that $\vol_d ( E_t )$ is smaller than   
$$    
\E \exp \Bigg( \alpha \Big( d s^2 - \sum_{j=1}^d (X_j-s)^2 \Big) \Bigg) =    
\Big(\E \exp \big( \alpha(2sX-X^2) \big)\Big)^d   
$$   
where $X$ is uniformly distributed in $[0,1]$ and $\alpha>0$ is arbitrary.   
This implies  
$$ \vol_d ( E_t ) \le \Big( \inf_{\alpha>0}  g(s,\alpha) \Big)^d $$   
where    
$$   
g(s,\alpha) = \int_0^1 \exp(\alpha(2sx-x^2)) \, \dint x.   
$$   
By continuity and the proof in \cite{DFM90} we find  a positive $t_0$, and    
for each $t\in[0,t_0]$, we find some positive $\alpha$ such that   
$$    
g(s,\alpha) < \frac{10}{11},    
$$   
where $2s=\tfrac12(1+t)$.    
Now \eqref{eq:_3} follows and the proof is completed.   
\end{proof}   
   
\section{$L_p$ approximation}    
   
The $L_p$ approximation problem is defined by  
$$  
{\rm APP}_d\,:\ F_d \to L_p([0,1]^d) \ \ \ \mbox{with}\ \ \  
{\rm APP}_d(f)=f  
$$  
for $F_d\in\{F_d^{\,\rm mon},F_d^{\,\rm con}\}$ and the standard  
$L_p([0,1]^d)$ space.   
  
The algorithms $A_n$ are now given by~\eqref{form} with  
$\phi_n:\reals^n\to L_p([0,1]^d)$. The $n$th minimal error for the  
$L_p$ approximation problem in the  
worst case setting is defined by  
$$  
e_n^{\rm app}(F_d)=\inf_{A_n}\ \sup_{f\in F_d}\|{\rm  
  APP}_d(f)-A_n(f)\|_{L_p([0,1]^d)}.  
$$  
For $n=0$, the initial error is again $\tfrac12$. The information  
complexity is now  
$$  
n^{\rm app}(F_d,\e)=\min\{n\,|\ \ e_n^{\rm app}(F_d)\le\e\}.  
$$  
  
Note that lower bounds for integration    
also hold for $L_p$ approximation.  
Indeed, take an arbitrary algorithm $A_n$ for the $L_p$ approximation  
problem, and let  
$$  
A_n^{\rm int}(f)=\int_{[0,1]^d}A_n(f)(x)\,{\rm d}x.  
$$  
Then $A_n^{\rm int}$ approximates the integral of $f$ and we have   
$$  
{\rm INT}_d(f)-A_n^{\rm int}(f)=  
\int_{[0,1]^d}\big(f(x)-A_n(f)(x)\big)\,{\rm d}x. 
$$  
This yields  
$$  
\left|{\rm INT}_d(f)-A_n^{\rm int}(f)\right|\le    
\int_{[0,1]^d}\left|f(x)-A_n(f)(x)\right|\,{\rm d}x\le  
\left(\int_{[0,1]^d}\left|f(x)-A_n(f)(x)\right|^p\,{\rm d}x\right)^{1/p}.  
$$   
Since this holds for all algorithms $A_n$, we have   
$$  
e_n^{\rm int}(F_d)\le e_n^{\rm app}(F_d)\ \ \ \mbox{and}\ \ \  
n^{\rm int}(F_d,\e)\le n^{\rm app}(F_d,\e),  
$$  
as claimed. In particular, the curse of dimensionality also holds    
for the $L_p$ approximation problem for both classes $F_d^{\,\rm mon}$  
and $F_d^{\,\rm con}$.   
   
   

\end{document}